\documentclass{notices}
\usepackage{amsfonts,amssymb,amsmath,amscd}
\usepackage{cancel}
\usepackage{graphicx}
\usepackage[dvipsnames]{xcolor}

\makeatletter
\g@addto@macro\normalsize{%
  \setlength\abovedisplayskip{4pt}
  \setlength\belowdisplayskip{4pt}
  \setlength\abovedisplayshortskip{0pt}
  \setlength\belowdisplayshortskip{4pt}
}

\renewcommand{\epsilon}{\varepsilon}
\renewcommand{\phi}{\varphi}
\renewcommand{\le}{\leqslant}

\title{Strategy flexibility in university mathematics}

\author{Peter H\"ast\"o
\affil{Peter H\"ast\"o is profesor of mathematics at the University of Helsinki, Finland. His email is peter.hasto@helsinki.fi}
\and
Jon R. Star
\affil{Jon R. Star is professor of education at Harvard University, USA. His email is jon\_star@gse.harvard.edu}
}

\begin{document}

\maketitle

Our colleague Miguel Ab\'anades asked Spanish first year university students to calculate 
$\frac d{dx} \log(x^2)$. Of 40 students, 39 solved the task using the chain rule: 
\[
\frac d{dx} \log(x^2) 
= 
\frac1{x^2} \frac d{dx} x^2 = \frac1{x^2} \cdot 2x = \frac 2x\,.
\] 
Only one student started by manipulating the logarithm: 
\[
\frac d{dx} \log(x^2)
=
\frac d{dx} (2\log x) = \frac 2x\,.
\] 
Both approaches, or solution strategies, are correct, but many mathematicians might consider the latter more ``clever''. 
The ability and propensity to use solutions tailored to features of a problem is what we call \textit{strategy flexibility}. 

This phenomenon has been observed over a long time (e.g.\ Wertheimer, 1945, see 
Section~\ref{sect:nature} below). In the past year two research literature reviews of the field have 
appeared (\cite{HonEtAl23} and \cite{Ver24}) each with over 70 references. We refer 
the reader to these surveys for a comprehensive overview of the research and 
comparisons of different conceptualizations of the phenomenon. As is usual, 
most of the research has been conducted at the elementary and middle school levels, 
even though the impact of strategy flexibility is perhaps even greater at higher levels. 

In this expository article, we make the case strategy flexibility is relevant 
for university mathematics and its teaching. 
We include concrete mathematical examples throughout the text which hopefully will 
allow the readers to connect the presented ideas with their own experience and 
perhaps draw some inspiration for their own teaching development journey.

We first present additional manifestations of flexibility or lack thereof 
as well as most research from the university level (Section~\ref{sect:examples}).
We also extrapolate from flexibility research at other 
levels and draw on our personal experience. 
Specifically, we describe mathematics education research related 
to flexibility (Section~\ref{sect:nature}) 
and provide some ideas for practical steps that can be taken 
in university classes (Section~\ref{sect:teaching}).
We argue that university mathematics teachers are 
well-positioned to cultivate the skills and attitudes necessary for 
strategy flexibility.
Finally, challenges and caveats to teaching flexibility are presented. 

We drop the word ``strategy'' from ``strategy flexibility'', since this is the only type of 
flexibility that we consider; see \cites{HonEtAl23, Ver24} for other  types flexibility in 
representations. We also note that some authors use the term ``adaptive'' 
instead of ``flexible''. 


\section{Examples of flexibility}\label{sect:examples}

Those who teach mathematics have undoubtedly 
observed the flip side of the coin: weaker students are prone to select unexpectedly 
inefficient solution strategies. As an example 
from this year's Finnish high-school leaving exam, far too many students 
used the quotient rule as follows
\begin{align*}
&\frac d{dx} \frac{x^3-15x^2+50x}{1000}
= \\
&\quad\frac{(3x^2-30x+50)\cdot 1000 - (x^3-15x^2+50x)\cdot 0}{1000^2}\,.
\end{align*}
Unsurprisingly, this approach exposed the students to a host of slip-up opportunities 
that are absent from a more efficient approach.

The opportunity to be clever --– to choose one's approach wisely --- occurs not only 
at university but at every level of mathematics. 
In arithmetic, one can calculate $5002-3997$ with a standard 
``borrowing''-approach or by noticing that adding $5$ to $3997$ gives 
$4002$, after which another $1000$ is needed to reach $5002$. 
Already back in 1945, Wertheimer lamented that students 
calculated
\[
\frac{274+274+274+274+274}5
\]
by first adding the five numbers in the numerator and then performing a division.
In elementary algebra, the equation $3(x-4)+4(x-4)=14$ can be solved by the 
standard algorithm which involves distributing both parentheses
or by a situationally appropriate strategy which starts with 
combining like terms:
\begin{align*}
3(x-4)+4(x-4)&=14 \\
7(x-4)&=14 \\
x-4&=2 \\
x&=6.
\end{align*}
University mathematics contains both similar opportunities for flexibility 
as well as more profound ones. 

In university calculus courses, we can use tasks like 
\[
\frac d{dx} \log(x^2) 
\qquad\text{and}\qquad
\frac d{dz} \big(z(z^3+z)^{-1}\big),
\]
from the first paragraph (above) and Maciejewski--Star's study \cite{MacS16} in the US, respectively. 
In both cases, a simplification before differentiation leads to a more efficient 
solution. 
This is analogous to the linear equation mentioned above. 
Schoenfeld (1985) makes the opposite observations about his students' approaches to integration tasks such as $\int \frac x{x^2-9}\,dx$: an initial algebraic manipulation 
(partial fractions) leads to much more work. 
Beyond calculus, the inefficient solution strategy may be so technically 
demanding that students are not able to complete it at all: 
a quantitative difference in performance becomes a qualitative one.

Broley and Hardy \cite{BroH22} studied tasks such as showing that the set
$\{\frac {n^p}{n+1} : n\in \mathbb N\}$ is unbounded when $p>1$ in an analysis course in Canada. 
They found that students mostly fell back on calculus routines and had difficulty 
connecting solutions to concepts from analysis. 
However, they suggested that flexibility in calculus could improve students' ability to make 
connections between calculus and analysis. 

In a linear algebra course at a New Zealand university, Kontorovich \cite{Kon20} 
similarly connected flexibility 
to the use of general properties over matrix calculation-based approaches. 
Students were for instance asked to find  a basis of the space of vectors $b\in \mathbb R^3$ 
for which the equation $Ax=b$ has a solution, where 
\[
A=\Bigg[\begin{array}{cc}1&2\\4&1\\0&1\end{array}\Bigg].
\]
In this case the most efficient approach was based on observing that the column vectors of 
$A$ are independent since their coordinates are not proportional 
(this was a lemma proved in the course). 
A computational approach was to use Gaussian elimination; this involved many more steps, which 
were completed in more or less efficient ways by different students. 


\section{The nature of flexibility}\label{sect:nature}

To recap, flexibility refers to having access to a variety of solution strategies and the ability and inclination to choose a 
\textit{situationally appropriate} one that is well-suited to the task at hand. 

Traditionally, mathematical proficiency has been framed with the constructs of procedural and conceptual knowledge. 
Procedural knowledge is often seen merely as the ability to execute memorized algorithms (reflecting superficial understanding) whereas conceptual knowledge is seen as indicating  deeper connections between mathematical objects. 
Star \cite{Sta05} argued that the procedural--conceptual dimension should be distinguished from the superficial--deep dimension. 
He proposed that flexibility is a form of deep procedural knowledge where one needs to understand the limitations and strengths of procedures in relation to specific tasks.

\begin{figure*}[htb]
\centering
\includegraphics[width=6.2in]{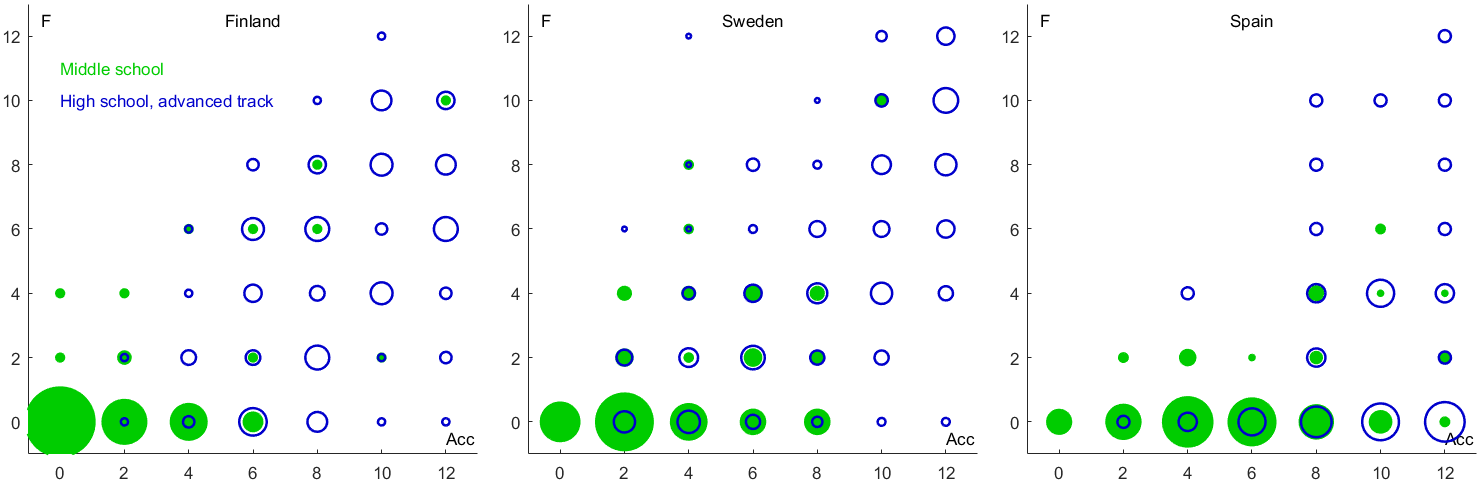}
\caption{Accuracy (``Acc'', horizontal) vs Flexibility (``F'', vertical) in Finland (left), Sweden (middle) and Spain (right) in middle school (solid green) and high school (blue outline). 
The area of the circle is proportional to the percentage of students with the given accuracy and flexibility scores. This previously unpublished figure is based on data from \cite{StaTJH22}.}
\label{fig:circles}
\end{figure*}

In investigations of elementary and secondary students' mathematics learning, there has been long-standing interest in this aspect of mathematical proficiency under a variety of labels, including productive thinking (Wertheimer), 
relational understanding (Skemp), 
adaptive expertise (Hatano), and 
strategic competence (Kilpatrick). 
Conceptualizing this competency as flexibility has gained traction in the field in the past 20 years (see \cites{HonEtAl23, Ver24} for recent reviews), perhaps due to its more precise definition, the development and use of validated measures of flexibility among researchers, and the ease of applicability of this construct to secondary as well as tertiary mathematics. 

One such instrument is the tri-phase test developed by Xu, Star and colleagues.
We used it in a fairly large study of middle- and high-school students' flexibility and accuracy in linear equation solving in three countries \cite{StaTJH22}. Unsurprisingly, we found 
students with low accuracy and low flexibility and others with high accuracy and high flexibility. Furthermore, many had high accuracy and low flexibility but virtually no one had high flexibility and low accuracy. This results in the lower-triangular pattern in Figure~\ref{fig:circles}. 
The fact that there were highly accurate students with all levels of flexibility indicates 
that flexibility is not merely a matter of having greater procedural skills or fluency; rather, flexibility involves a different kind of (deeper) procedural knowledge. 

The deep knowledge is needed since deciding on the 
situationally appropriate strategy requires paying attention to the special structure of a  problem and assessing the importance of these features in relation to the problem-solving objective \cite{Sta05}. 
While one might think that most problems have generic rather than special structure, 
experience suggests that problems with special structures are far more common in mathematical research and applications
than their ``statistical probability'' would suggest (e.g.\ if you generate polynomials with random coefficients, then clever shortcuts will almost never apply). 
Indeed, things lining up, canceling and simplifying often indicates to a mathematician that they are on the right track to solving the problem. 

The ability to view objects and expressions as the combination of parts 
has been called \textit{structure sense} by researchers. This concept also appears 
in the Common Core High School objective ``Seeing Structure in Expressions''. 
To summarize the previous paragraph: flexibility requires structure sense and 
flexibility also allows students to capitalize on their structure sense by choosing 
the most suitable strategy to the structure at hand. 
Most research on flexibility has dealt with primary and secondary school. Yet structure sense often comes into play far beyond arithmetic and linear equations.
Indeed, much of calculus involves pattern-matching expressions with formulas and carrying out correct algebraic manipulations. However, structure sense carries even further. For instance, in real analysis Hölder's inequality
\[
\int f(x)g(x)\, dx\le \|f\|_p\|g\|_{p'}\,, \quad\frac1p+\frac1{p'}=1, 
\]
is simple enough for students to understand, but applying it in the proof of Minkowski's inequality requires interpreting the integrand of $\int |u|^p\, dx$ as $|u|\cdot |u|^{p-1}$ and identifying the first term with $f$ and the second term with $g$ in Hölder's inequality. This kind of application is much more challenging, especially for students who are less adept at spotting patterns in algebraic expressions.

Flexibility involves not only skills and knowledge; it is also a matter of disposition. 
Core to flexibility is striving for the most appropriate or best solution in a given situation. Some researchers have sought to quantify ``best'' in terms of the speed and accuracy with which the strategy can be applied by a given individual to a given problem. A good summary of this point-of-view is 
given by Hickendorff \cite{Hic22}. 
While a situationally appropriate strategy often involves fewer steps and can be executed quicker, we consider this a side-effect rather than the goal of flexibility. 
It is more important to contemplate a variety of approaches and weigh their pros and cons, rather than going with the first idea that springs to mind.

Let us speculate on more wide-ranging consequences of such a mindset of optimization 
beyond tasks in a mathematics lesson. Expanding on the point that it impacts all 
levels of mathematics, even reading research papers one sometimes wonders if this is just the first proof that the authors managed to produce and if additional effort could have led to  
more streamlined and elegant proofs. 
In programming, there is a world of difference between a brute-force solution and an elegant one. 
Choosing  data-structures unwisely will stifle opportunities to evolve one's program with changing demands. 
More generally still, one can assess and improve most processes ranging from engineering to administration in terms of efficiency if one spends time and effort on optimization.

The critical reader may wonder whether all this ``flexibility'' is actually a manifestation 
of the same trait. 
We see a conceptual relationship between these different examples of 
flexibility and at least personally utilize the 
same flexible thinking also outside procedural mathematics.
However, research suggests that transfer of skills between contexts is usually difficult and 
the jury is still out on the level of transfer of (acquired) flexibility to other 
mathematical areas and other subjects. 
What the research does demonstrate is that it is possible to impact 
students' flexibility positively by suitable adjustments to teaching 
practices. 

%

\section{Teaching for flexibility}\label{sect:teaching}

Much has been written about mathematical problem solving (in the sense of P\'olya), 
including at the university level (e.g.\ Schoenfeld, 1985);  
teaching problem solving skills is likely to also support flexibility. 
However, despite decades of effort it is fair to say that teaching problem solving 
has not had a large-scale impact, due in part to disagreement on what ``problem solving'' is. 
Compared to the daunting task of improving mathematical 
problem solving, aiming for flexibility offers several advantages. 
Both a traditional routine-oriented approach and flexibility involve a focus on
procedural knowledge. However, adding an element of flexible procedure use changes students' perceptions of the mathematical task, from merely parroting back what a teacher has demonstrated to a wider, yet still limited, set of creative options. 
It is the limited range of creativity which sets flexibility apart from problem solving. 
Thus the reorientation required of the teacher and the students is smaller in 
flexibility than in mathematical problem solving.

Star and Newton \cite{StaN09} studied mathematics experts'
flexibility. Participants recounted that looking for a situationally appropriate strategy in a mathematical task is intellectually challenging and interesting. 
Furthermore, they reported not having been taught to search for optimal strategies, but rather considered doing so a personal desire that they brought to problem solving. 
While these mathematics experts had been able to figure out flexible strategy use on their own, 
earlier mentioned studies have shown that this is not the case for most 
students.

In reference to Figure~\ref{fig:circles}, we mentioned that there were virtually no 
students with high flexibility and low accuracy in our study \cite{StaTJH22}. 
The figure indicates that the Spanish schools in the sample placed more emphasis on 
accuracy than on flexibility compared to the Nordic countries. This is evidence that 
it is possible to impact flexibility by teaching and curriculum choices.
We are not aware of similar cross-national comparisons of flexibility at the 
university-level. However, Ernvall-Hytönen, Hästö, Krzywacki and Parikka \cite{ErnHKP22} 
found that Finnish university students showed fair levels of flexibility whereas 
Shaw, Pogossian and Ramirez \cite{ShaPR20} found that students in a 
US university showed less flexibility than they expected. It should be mentioned that 
these studies measured flexibility with very different instruments. 
Nevertheless, variability across contexts 
suggests that flexibility is not (only) 
a consequence of being mathematically gifted, but rather 
depends on the education system and can be cultivated through teaching. 
It is important not to view flexibility as the ability to think of unusual, clever solutions 
out of nowhere (``fixed mindset''), since such a view can be dispiriting and counter-productive. 
Rather, flexibility is a combination of skills and attitudes that can be achieved through 
systematic work and supported through instruction (``growth mindset'').

Star and Seifert \cite{StaS06} carried out an intervention which illustrates 
a fairly straight-forward yet generalizable teaching approach. Students were introduced 
to equation solving by teaching them what operations are permitted on equations 
(add constant to both sides, combine like terms, etc.)\ but were not told to follow a specific 
sequence of steps or rules (such as the standard algorithm). 
Furthermore, students were instructed to solve each problem with at least two 
different strategies. 
This focused the attention on the solution process instead of the answer 
and showed the students that there are multiple viable solution strategies. 
As a result, the intervention participants exhibited more varied 
solution strategies while also achieving the same proficiency on the standard 
algorithm when compared to the control group. 
In a calculus context, a corresponding shift could be to tackle multiple 
differentiation rules simultaneously with tasks that match several rules at the 
same time. An additional layer is that calculus also involves 
algebraic manipulations that may invite flexibility; the examples in the introduction utilize this. 
A one-session intervention along these lines at a Canadian university was studied by 
Maciejewski and Star \cite{MacS16}.

The students in the interventions \cites{MacS16, StaS06}
had more freedom to work things out than in traditional teaching, 
but the number of options offered to them was not overwhelming. 
Furthermore, instead of instilling students with 
an ``always distribute first''-mantra, this approach to equation solving encourages them 
to pause and think for a while before plunging into solving a problem. 
Thinking before acting should help reduce the frequency of 
strangely inefficient solutions that surface far too often, 
such as solving $(x+1)(2x-3)=0$ 
by distributing and using the quadratic formula on $2x^2-x-3=0$.
Scaling back on the ``always follow this set of steps''-approach also 
enables different task types where students reverse the steps, 
start in the middle, analyze given solutions, etc. 
For instance, students could start with an answer and then 
manipulate it to construct an equation for the rest of the class to solve; 
the process could look like
\begin{align*}
x&=6 \\
3x&=18 \\
3x-4&=14 \\
\frac12(3x-4)&=7.
\end{align*}
This task highlights the invertibility of the steps in equation solving. 
An analogous task in calculus could elucidate the relationship between 
differentiating and integrating. To avoid overly trivial tasks, students should be encouraged to aim for ``sneaky'' problems that will give their classmates a run for their money. 

Although flexibility is sometimes framed as using the optimal choice among a set of known 
strategies, the intervention \cite{StaS06} also shows that it is possible for students to construct 
their own strategies. As long as the work is done in a suitably restricted domain, the complexity and cognitive load can be kept at reasonable levels. 
The fewer the restrictions, the closer the situation is to problem solving. 

We emphasize that flexibility is more than just generating multiple strategies: the strategies should also be compared and evaluated using mathematical criteria. The teacher is instrumental in negotiating the socio-mathematical norms 
with the class regarding the appropriate criteria to use in the comparison as well as making sure that every student's contribution is appreciated while still discussing different proposed strategies critically. 
Comparison between solutions can be carried out with any task for student-generated 
solutions.
%
Furthermore, selecting the structure of the formulas and numbers appropriately can encourage applying an opportune, non-generic strategy. For instance, the integral 
\[
\int_1^2 x^2-2x+1\, dx
\]
invites the change of variables $z:=x-1$, even though a direct calculation is also quite feasible. 

On the other hand, it is also possible to take an explicit flexibility focus in a task. 
For example, students can be asked to provide two (or more) different solutions to a given 
problem. The comparison of the solutions can be part of the task or organized by the teacher 
as a classroom discussion. 
Supporting flexibility by digital tools is another interesting, yet underdeveloped 
area \cite{BreMVH19}.

\begin{figure*}
\centering
\textbf{Problem}: Calculate $\displaystyle\frac d{dz} \big(z(z^3+z)^{-1}\big)$.
\medskip

\begin{tabular}{p{3.5cm}p{3.8cm}@{\hskip 10pt}||@{\hskip 12pt}p{3cm}p{3.5cm}}
\textbf{Avery's solution}: &&& \textbf{Blake's solution}:  \\
I rewrite the expression as a quotient &$\displaystyle\frac d{dz} \frac z{z^3+z}$ 
& $\displaystyle z^3+z = z(z^2+1)$ & I factorize the denominator\\
I use the quotient rule & $\displaystyle=\frac{(z^3+z)-z(3z^2+1)}{(z^3+z)^2}$ 
& $\displaystyle\frac d{dz} (z^2+1)^{-1}$ & I cancel the $z$ \\
I simplify the numerator and 
factorize the denominator & $\displaystyle=\frac{-2z^3}{(z^2+1)^2z^2}$
& $\displaystyle=-(z^2+1)^{-2} \frac d{dz}z^2$ & I use the power rule and the chain rule \\
I cancel the $z^2$ & $\displaystyle=-\frac{2z}{(z^2+1)^2}$\,.
& $\displaystyle=-2z(z^2+1)^{-2} $ & I use the power rule again. \\
\end{tabular}
\begin{itemize}
\centering
\setlength{\itemsep}{0pt}
\item
Highlight the similarities and differences between the solutions.
\item
What are the advantages and disadvantages of each strategy?
\end{itemize}
\caption{A worked example pair}
\label{fig:WEP}
\end{figure*}

A task format specifically designed to support flexibility is the worked example pair (WEP). 
A WEP presents two hypothetical students' solutions to a problem (see Figure~\ref{fig:WEP}). In the simplest case 
one solution represents the generic approach and the other is situationally appropriate. 
The WEP also contains some prompts such as:
\begin{itemize}
\setlength{\itemsep}{0pt}
\item
Explain what is done in the solutions. 
\item
Explain why this is a valid solution to the problem.
\item
Highlight the similarities and differences between the solutions.
\item
What are the advantages and disadvantages of each solution?
\item
Which solution is better?
\item
Can you solve the problem in yet another way?
\item
Which solution is correct?
\end{itemize}
As the last question hints at, one subtype of WEPs involves an error in 
one or both of the solutions. Such intentional errors may help shift students' attention from 
the answer to the argumentation and lessen their anxiety about presenting in class 
a solution with an error \cite{PalH18}. Furthermore, debugging a solution is an 
important learning opportunity in itself. 

Students' attitudes are a potential challenge to a focus on flexibility. Prior experience from the education system may have led to the belief that the most, or only, important thing in mathematics is to obtain the correct answer, and so discussions about multiple solutions are a waste of time. 
Students may believe that it is the teacher's job to explain the correct strategy for a problem-collection and their job is to apply it as instructed. Such beliefs hinder teaching higher levels of thinking, 
including flexibility. 

In the short-term, an answer-centered orientation may lead to seemingly good results. 
When students are drilled on a set of tasks for a few weeks, 
most will display adequate performance on a test with similar tasks at the end of the training period. However, transfer of the knowledge to different situations, 
application and retention to a later time are likely to be quite limited (Boaler, 1993). 
In a middle-school intervention for flexible equation solving in Finland 
many teachers suggested the initial slow progress was more than compensated 
for by not needing to start reviewing equations from scratch in subsequent years.
Good short-term test results arise from students being relieved of many processes involved in solving the problem, such as identifying the problem type and relevant solution methods 
or combining techniques and concepts from several areas. 
But these are important habits of mind that the students should be learning. In the long-term, this leads to worse structure sense and weaker meta-cognitive skills such as the ability to assess the viability of a strategy in a situation. 

One source of such unhelpful beliefs is assessment. In some education systems and institutions, 
teachers assess students' work by only scoring the correctness of the final answer. 
Scoring only the answer implicitly communicates to students that the argument is not important. In secondary and tertiary education, points are often allocated also to the argument. 
However, in mathematics usually any correct argument gives full points, irrespective of its other qualities such as elegance, parsimony or readability. 
In particular, the immediate pay-off from flexibility may be quite low, resulting only from occasionally saving a little time on calculations. Thus students may need additional incentives to strive 
for flexibility and the concomitant small advantages that accumulate over time. 
For this reason, Maciejewski \cite{Mac22} suggested including some tasks 
in calculus exams which directly 
relate to flexibility, such as asking for multiple solutions to the same problem or using 
a WEP. More radically still, would you consider giving different points for 
Avery's and Blake's solution strategies in Figure~\ref{fig:WEP}?
This seems like a strange idea in mathematics, whereas in most other subjects an elegant presentation 
of an argument would be given more points than a cumbersome one.

%
%
%


\section{Concluding summary}\label{sect:summary}

We have argued that procedural knowledge can have more depth than it seems at first glance, specifically when it comes to choosing or constructing a solution strategy tailored to a specific problem. The ability and inclination to use situationally appropriate strategies is called strategy flexibility, and it has been an active area of mathematics education research in the past two decades with roots stretching much further back.

Research has shown that students across skill-levels can learn flexibility with appropriate support and encouragement. Flexibility is furthered by shifting the focus from answers to solutions strategies and the ideas behind them. 
Most tasks can be used to support flexibility; however, there are also tasks specifically designed to promote flexibility, such as worked example pairs. 
Most research has focused on flexibility pre-university, but many of the tools and task-types can be applied 
in a university context. University mathematics teachers are well-positioned 
to appreciate the value of flexibility and convey it to their students, which may be the 
most important factor for adopting a flexible mindset.

Flexibility provides a framework for teaching mathematics with a balance of conceptual and procedural focus. It allows teachers and students to expand their working habits gradually toward a more problem-solving oriented style. Teaching for strategy flexibility is a small step for the teacher, but it may just be a small leap for students who realize that mathematics involves decisions and multiple correct solution strategies with different strengths.


\section*{Acknowledgment}

We thank Anne-Maria Ernvall-Hytönen for comprehensive feedback on this paper and
Heidi Krzywacki, Päivi Portaankorva-Koivisto as well as the referees for comments. 


\bibliographystyle{amsplain}

%
%
%
\end{document}